\documentclass[10pt,twoside]{amsart}
\usepackage{amsmath,amsfonts,amssymb,amsxtra,setspace,xspace,graphicx,lmodern,psfrag,epsfig,color,latexsym}
\usepackage[T1]{fontenc}
\usepackage[colorlinks=true]{hyperref}\hypersetup{urlcolor=blue, citecolor=red}

\newcommand{\bbZ}{\mathbb{Z}}

\usepackage[textsize=small,backgroundcolor=red,linecolor=red!30,textwidth=2cm]{todonotes}

\newcommand{\R}{\mathbb{R}}

\renewcommand{\vec}[1]{{\boldsymbol#1}}

\newcommand{\ie}{\textit{i.e.}\/, }
\newcommand{\eg}{\textit{e.g.}\/, }

\def\e{\epsilon}

\def\O{\Omega}

\def\e{\eta}

\def\ue{u_\e}

\def\uek{u_{\e_k}}

\def\uo{u_0(\vec{x},\vec{y})}

\def\ve{v_\e}

\def\wdto{\buildrel\hbox{${\bf R}$}\over
\rightharpoonup \!\!\!\!\! \rightharpoonup}

\def\dto{\buildrel\hbox{${\bf R}$}\over
\to \!\!\!\!\! \to}

\def\tol2{\buildrel\hbox{$L^2$}\over\longrightarrow}

\def\IOm{\int_{\Omega}}

\def\IOAY{\int\!\!\int_{\Omega\times Y^m}}

\def\limf{\displaystyle{\liminf_{\eta\to 0}}}
\def\lims{\displaystyle{\limsup_{\eta\to 0}}}

\def\xxe{\left(\vec{x},\frac{{\bf R}\vec{x}}{\e}\right)}
\def\xy{(\vec{x},\vec{y})}

\newtheorem{definition}{Definition}[section]
\newtheorem{theorem}{Theorem}[section]
\newtheorem{proposition}{Proposition}[section]
\newtheorem{lemma}{Lemma}[section]
\newtheorem{remark}{Remark}[section]
\newtheorem{corollary}{Corollary}[section]

\begin{document}
\let\WriteBookmarks\relax
\def\floatpagepagefraction{1}
\def\textpagefraction{.001}



\title[Two-scale convergence for quasiperiodic monotone operators]{Two-scale cut-and-projection convergence for quasiperiodic monotone operators}

\bibliographystyle{alpha}

\author{Niklas Wellander$^{1}$}
\address{$^{1}$Swedish Defence Research Agency (FOI), FOI P.O. Box 1165 SE-581 11, Link{\"o}ping, Sweden}

\author{Sebastien Guenneau$^{2}$}
\address{$^{2}$UMI 2004 Abraham de Moivre-CNRS, Imperial College London, SW7 2AZ London, UK}

\author{Elena Cherkaev$^{3}$}
\address{$^{3}$University of Utah, Department of Mathematics, 155 South 1400 East, JWB 233, Salt Lake City, USA}

\begin{abstract}
Averaging certain class of quasiperiodic monotone operators can be simplified to the periodic homogenization setting by mapping the original quasiperiodic structure onto a periodic structure in a higher dimensional space using cut-and projection method. We characterize cut-and-projection convergence limit of the nonlinear monotone partial differential operator $-\mathrm{div} \; \sigma\left(\vec{x},\frac{{\bf R} \vec{x}}{\e}, \nabla  u_\e\right)$ for a bounded sequence $u_\e$ in
$W^{1,p}_0(\Omega)$,  where $1<p < \infty$, $\Omega$ is a bounded open subset in $\R^n$ with Lipschitz boundary. We identify the homogenized problem with a local equation defined on the hyperplane in the higher-dimensional space. A new corrector result is established.  
\end{abstract}




\maketitle

\section{Introduction}\label{intro}
Nonlinear physical phenomena are ubiquitous in modern electronic devices. A few examples are  current surge protectors made of varistor ceramics, solid state amplifiers and integrated circuits. This is one motivation to develop mathematical tools that can be used to analyse the effective properties of polycrystalline quasiperiodic semiconductors. 

In \cite{braides2009homogenization}, it is shown that integral energies $F_\eta$ where the spatial dependence follows the geometry of a Penrose tiling, or more general quasicrystalline geometries, can be homogenized. More precisely,
\begin{equation}
F_\eta(\vec{x})=\int_\Omega f\left(\frac{\vec{x}}{\eta},\nabla  u(\vec{x})\right) \, d\vec{x} \; , \; u \in W^{1,p}(\Omega)
\label{gamma1}
\end{equation}
where $\Omega$ is an open subset of $\R^2$, and $f$ depends on $\vec{x}$ through the shape and the orientation of the cell containing
$\vec{x}$ in an a-periodic tiling of the space,
$\Gamma$-converge on $W^{1,p}(\Omega)$ with respect to the $L^p$ convergence to the functional
\begin{equation}
F_0(\vec{x})=\liminf_{T\to\infty}\left\{\frac{1}{T^2}\int_{(0,T)^2} f\left(\vec{y},\nabla  v(\vec{y})+\xi\right) \, d\vec{y} \; , \; v \in W^{1,p}_0((0,T)^2)\right \}
\label{gamma2}
\end{equation}
This general homogenization result was shown using that $f$ is Besicovitch almost periodic in $\vec{y}$ and thus a previous result on Besicovitch almost periodic functionals \cite{braides1986homogenization} could be applied. Homogenization of interfacial energies on Penrose lattices making use of $\Gamma$-convergence for similar functionals to (\ref{gamma1}) but with the volume integral replaced by a surface integral has been also addressed in \cite{braides2012interfacial}.

$\Gamma$-convergence is a very powerful tool in homogenization theory \cite{braides2002gamma}, but two-scale convergence \cite{Nguetseng1989,Allaire1992} can more easily identify homogenized equations in the periodic setting. A similar tool is the periodic unfolding approach \cite{Cioranescu+etal2008}, in which one first maps the original sequence of functions to a sequence that is defined on $\R^n \times ]0,1[^n$, and then takes the usual weak limit in suitable function spaces, using this extended domain. This is similar to the two-scale Fourier transform approach proposed in \cite{Wellander2009}.

Due to its simplicity, one might wish to apply two-scale homogenization (or the periodic unfolding and Fourier transform approaches) to quasiperiodic materials or mixtures of materials with rational and irrational periodicity, \eg see \cite{Braides1991,Diaz+Gayte2002} for a setting in an almost periodic regime \cite{Blanc+etal2015a, Blanc+etal2015b} for some recent work on quasiperiodic multiscale homogenization setting.

This has been proposed in \cite{guenneau2001homogenisation,Bouchitte+etal2010} wherein two-scale convergence is applied to the quasi-periodic setting making use of the cut-and-projection method.
Indeed, quasiperiodic materials can be often described by periodic structures in higher spatial dimensions cut and projected onto a  hyper-plane or a lower dimensional (physical) space, typically $\R^2$ (such as for Penrose tilings) and $\R^3$, as proposed by the physicists Duneau and Katz \cite{Duneau+Katz1985}. This makes it possible to use standard periodic homogenization tools such as two-scale convergence, to homogenize quasiperiodic materials.

To do that one has to complement existing tools with the cut-and-projection operator, this was done in \cite{Bouchitte+etal2010} in the framework of $W^{1,2}$, making use of Fourier representation for two-scale limits of gradients. Importantly, this was revisited in \cite{Wellander+etal2017,Wellander+etal2019}.

However, two-scale convergence can  also be applied   to nonlinear operators \cite{Allaire1992}. 
This is a generalization of the usual weak 
convergence in Lebesgue spaces $L^p$, $1<p <  \infty$, in which one uses oscillating test functions to capture oscillations on the same scale as the test functions in the sequence of functions that are investigated. As a consequence one obtains limit functions that are defined on the product space $\R^n \times ]0,1[^n$.

 In this paper we extend the two-scale cut-and-projection convergence method to Sobolev spaces $W^{1,p}$, $1 < p < \infty$.  We build upon \cite{Wellander+etal2019} to characterize the limits of nonlinear partial differential operators in this setting. We  illustrate the method on a nonlinear electrostatic problem that was previously homogenized using the tool of G-convergence for a larger class of almost periodic functions in \cite{Braides+etal1992}. We finally establish a corrector result for the gradients. We note that Kozlov established in \cite{kozlov1979averaging} a corrector result that requires almost periodic coefficients, which are the restriction of sufficiently smooth periodic functions of greater
number of variables, when the problem is set in the whole space.

\subsection{Setup of a nonlinear electrostatic equation in a quasiperiodic structure} 

Throughout the paper, we consider a bounded domain $\Omega$ in $\R^n$ with Lipschitz boundary. We study the electrostatic equation, which is  applicable to model DC currents in  semiconductors, \eg   ZnO based varistor ceramics, 
 \begin{equation}  \label{eq:NL_static}
     \left\{
        \begin{aligned}
          -  \mathrm{div} \; \sigma\left( \vec{x}, \frac{{\bf R}\vec{x}}{\e}, \nabla u_\e(\vec{x}) \right)& = f (\vec{x}) \;, \qquad \vec{x} \in \Omega\\
        \left.    u_\e\right|_{\partial\Omega} & = 0
        \end{aligned}
    \right.
 \end{equation}
with $f \in W^{-1,q}(\Omega)$ and $p$ is a real constant $1<p <  \infty$, and $q$ is its dual exponent, \ie  $1/p+1/q=1$. We use standard notations for Lebesgue and Sobolev spaces.  The Euclidean norm and the scalar product in $\R^n$ are denoted by $\mid \cdot  \mid$ and $(\cdot ,\cdot )$, respectively. 
The current density is given by a non-linear map, $\sigma$, that satisfies  assumptions  i)-vii):
\begin{itemize}
    \item[i)] $\sigma(\vec{x}, \cdot, \xi)$ is $Y$-periodic in $\R^m$ and Lebesgue measurable for every $\vec{x}, \vec{\xi}\in\R^n$.
    
    \item[ii)] $\sigma(\cdot, \vec{y}, \xi)$ is continuous for almost every $\vec{y} \in \R^m$ and  every $\vec{\xi}\in\R^n$.
    
    \item[iii)] $\sigma(\vec{x},\vec{y}, \cdot)$ is continuous for almost every $\vec{x}\in\R^n$ and $\vec{y} \in \R^m$. 
    \item[iv)]   $ 0 \leq c  |\vec{\xi}|^p \leq  \left( \sigma(\vec{x},\vec{y},\vec{\xi}) , \vec{\xi} \right) \, ,   c>0 $     for almost every $\vec{x}\in\R^n$  and $\vec{y} \in \R^m$.
    
    \item[v)]  $\left(\sigma(\vec{x}, \vec{y}, \vec{\xi}_1)- \sigma(\vec{x}, \vec{y}, \vec{\xi}_2) , 
    \vec{\xi}_1-\vec{\xi}_2 \right)\geq c_1 \left|\vec{\xi}_1 - \vec{\xi}_2 \right|^p, \; c_1>0 $ for all $ \vec{\xi}_1,\vec{\xi}_2  \in \R^n$, and almost every $\vec{x}\in \Omega$  and $\vec{y}\in\R^m$.
    \item[vi)]  $|\sigma(\vec{x}, \vec{y}, \vec{\xi}) | \leq c_2 \left(1 + |\vec{\xi}|^{p-1}\right), \; c_2>0, \forall \vec{\xi} \in \R^n$, almost every $\vec{x}\in \Omega$  and $\vec{y}\in\R^m$.
    
    \item[vii)]  $\sigma(\cdot , \mathbf{R} \cdot/\e, \xi) $ is Lebesgue measurable, for all $\vec{\xi} \in \R^n$. 
\end{itemize}
Standard estimates yield  solutions that are uniformly bounded in $W^{1,p}_0(\Omega)$ with respect to $\e$.  
\begin{remark}
Assumption v) is needed for a corrector result. 
It can be replaced by a strict monotone assumption if only the homogenized equation is needed, \eg see \cite{Allaire1992}. 
\end{remark}

\subsection{Two-scale cut-and-projection convergence}
In this section, we recall some properties of two-scale convergence in $L^p(\Omega)$, $1<p < \infty$,   
$\Omega\subset\R^n$ \cite{Allaire1992}, and revisit the extension to the quasiperiodic setting, see \cite{Bouchitte+etal2010} when $p=2$. More precisely, we consider a real valued matrix ${\bf R}$ with $m$ rows and $n$ columns, and we would like to approximate an oscillating sequence $\{u_\e(\vec{x})\}_{\eta\in ]0,1[}$ by a sequence of two-scale functions $u_0\left({\vec{x}},\frac{{\bf R}\vec{x}}{\e}\right)$ where $u_0\left({\vec{x}},\cdot\right)$ is $Y^m$-periodic on $\R^m$, a 
higher dimensional space paved with periodic cells $Y^m={]0,1[}^m$. In what follows, we assume that ${\bf R}: \R^n \to \R^m$, $n < m$, fulfills the criterion
\begin{equation} \label{criterion}
{\bf R}^T \vec{k} \not = {\bf 0} \; , \; \forall \vec{k} \in
\bbZ^m \setminus \{{\bf 0}\}  
\end{equation} 
In fact, when defining a quasiperiodic structure through cut-and-projection, one notes that the matrix ${\bf R}$ is not uniquely defined (\eg an icosahedral phase using a mapping from $\R^6$ or $\R^{12}$ onto $\R^3$, or a Penrose tiling using a mapping from $\R^4$ or $\R^{5}$ onto $\R^2$, 
see \cite{Duneau+Katz1985,Janot1992}). However, if $g$ is a trigonometric polynomial, then the $f = g \circ {\bf R}$ admits the following
(uniquely defined) ergodic mean:
\begin{equation}
L(f) = \lim_{T\to + \infty} \frac{1}{{(2T)}^n}\int_{{]-T;T[}^n} f(\vec{x}) \, \mathrm{d}\vec{x} = \int_{Y^m} g(\vec{y}) \, \mathrm{d}\vec{y} = [g] 
\label{ergodic}
\end{equation}
 where $[g]$ denotes the mean of $g$ over the periodic cell $Y^m$ in $\R^m$. This is the case
 provided that ${\bf R}$ fulfills the criterion (\ref{criterion}), see \cite{Bouchitte+etal2010}.

This result suggests the following concept of two-scale convergence
attached to a matrix ${\bf R}$.

\begin{definition}[Distributional two-scale convergence]
We say that the sequence $(u_\e)$ in $L^p(\Omega)$,
$1<p < \infty$, 
two-scale converges in the distributional sense towards the
function $u_0\in L^p(\O\times Y^m)$ for a matrix ${\bf R}$,  
if for every $\varphi\in{\mathcal{D}}(\O;C^\infty_\sharp(Y^m))$:
  \begin{equation}\label{eq:distr2scale}
\lim _{\e\to 0} \IOm\ue (\vec{x}) \varphi\left(\vec{x},\frac{ {\bf
R}\vec{x}}{\e}\right) \,\mathrm{d}\vec{x} = \IOAY \uo \varphi(\vec{x},\vec{y})\,\mathrm{d}\vec{x}\mathrm{d}\vec{y}
  \end{equation}
  \label{2scaledistribution}
\end{definition}


\begin{definition}[Weak two-scale convergence]
We say that the sequence $(u_\e)$ in $L^p(\Omega)$,
two-scale converges weakly towards the function $u_0\in L^p(\O\times Y^m)$ for a matrix ${\bf R}$,  
if for every $\varphi\in L^q(\O,C_\sharp(Y^m))$, $1 < p < \infty$, $1/p + 1/q =1$,  \eqref{eq:distr2scale} holds.
\label{2scaleweak}
\end{definition}

We denote weak two-scale convergence for a matrix ${\bf R}$ with   $\uek \wdto u_0$.
The following result, which is a straightforward extension of a proof in \cite{Bouchitte+etal2010} to $L^p$ case corresponding to Corollary 1.15 in  \cite{Allaire1992}, ensures the existence of such two-scale limits when the sequence $(u_\e)$ is bounded in $L^p(\O)$ and ${\bf R}$
satisfies (\ref{criterion}). 

\begin{proposition}\label{prop:weaktwoscaleprop}
If ${\bf R}$ is a matrix satisfying (\ref{criterion}) and $(\ue)$ is a
bounded sequence in $L^p(\O)$, $1< p < \infty$, then there exist a vanishing
subsequence $\e_k$ and a limit $\uo \in L^p(\O\times Y^m)$
($Y^m$-periodic in $\vec{y}$) such that $\uek \wdto u_0$, as $\eta_k \to 0$.
\end{proposition}
A proof of Proposition~\ref{prop:weaktwoscaleprop} uses the same arguments as in the periodic case, \eg see  \cite{lukkassen2002two} and can be found in   \cite{ferreira2021homogenization}.

We will need to pass to the limit in integrals
$\displaystyle{\int_{\Omega} u_\e \; v_\e \, \mathrm{d}\vec{x}}$ where
$u_\e\wdto u_0$ and $v_\e\wdto v_0$. For this, we introduce 
the notion of strong two-scale (cut-and-projection) convergence for a matrix ${\bf R}$.
\begin{definition}[Strong two-scale convergence]\label{def:strong_two_scale}
A sequence $u_\e$ in $L^p(\O)$ is said to two-scale converge strongly, for a matrix ${\bf R}$, towards a limit $u_0$ in
$L^p_\sharp (\O\times Y^m)$, which we denote $u_\e\dto u_0$, if and
only if $\uek \wdto u_0$ and 
\begin{equation}
\displaystyle{ {\Vert u_\e(\vec{x})\Vert}_{L^p(\O)} \to{\Vert u_0
(\vec{x},\vec{y}) \Vert}_{L^p(\O\times Y^m)}}  
\end{equation}
\end{definition}

This definition expresses that the effective oscillations of the
sequence $(u_\e)$ are on the order of $\e$. Moreover,
these oscillations are fully identified by $u_0$.

We note the following result which is useful to establish a link between weak quasiperiodic convergence in definition \ref{2scaleweak} and strong $L^p$ convergence (see Corollary \ref{L2to2scale}). 
\begin{proposition}\label{L1convergence}
Let ${\bf R}$ be a linear map from $\R^n$ to $\R^m$ satisfying
\eqref{criterion}. Let $f(\vec{x},\vec{y})\in L^1(\O; C_\sharp(Y^m))$.
Then $f\xxe$ is a measurable function on $\Omega$ such that
\begin{equation}\label{toto1}
\displaystyle{{\left\| f\xxe \right\|}_{L^1(\O)}
\leq
{\left\| f(\vec{x},\vec{y}) \right\|}_{L^1(\O; C_\sharp(Y^m))}
:=\int_{\Omega} \sup_{\vec{y}\in Y^m} \mid f(\vec{x},\vec{y}) \mid \, d\vec{x}}
\end{equation}
and
\begin{equation}\label{bgz2010}
\lim _{\e\to 0} \IOm f\left(\vec{x}, \frac{{\bf R}\vec{x}}{\e}\right) \,\mathrm{d}\vec{x} = \IOAY f(\vec{x},\vec{y})\,\mathrm{d}\vec{x}\mathrm{d}\vec{y}
  \end{equation}
\end{proposition}

{\bf Proof}  
Measurability of $f\xxe$ follows from Theorem 1 in \cite{lukkassen2002two} that ensures $f$ is of Carathedory type.
Inequality (\ref{toto1}) is obvious.  
Finally, (\ref{bgz2010}) follows from Lemma 2.4 in \cite{Bouchitte+etal2010}. 
\hfill \boxed{}

\vspace{5mm}
We deduce two corollaries.  First, we have a result that is useful to establish a corrector type result (see proposition \ref{allaireLp}). 

\begin{corollary}\label{Lpconvergence}
Let ${\bf R}$ be a linear map from $\R^n$ to $\R^m$ satisfying
\eqref{criterion}. Let $\phi(\vec{x},\vec{y})\in L^p(\O; C_\sharp(Y^m))$.
Then
\begin{equation}\label{toto2}
\lim _{\e\to 0} \IOm {\left| \phi\left(\vec{x}, \frac{{\bf R}\vec{x}}{\e}\right)\right|}^p \,\mathrm{d}\vec{x} = \IOAY {\left|\phi(\vec{x},\vec{y})\right|}^p\,\mathrm{d}\vec{x} \mathrm{d}\vec{y}
  \end{equation}
\end{corollary}
{\bf Proof}
We first consider the case when $\phi$ can
be expressed as $\phi\xy=\tau(\vec{x})\beta(\vec{y})$, where $\tau(\vec{x})\in L^\infty(\Omega)$ and $\beta(\vec{y})\in C_\sharp(Y^m)$. Since
$\beta^p\in C_\sharp(Y^m)$, we deduce from proposition \ref{L1convergence} that $\beta^p\left(\frac{{\bf R}\vec{x}}{\eta}\right)$ converges towards its mean $[\beta^p]$ weakly in $L^1(\Omega)$.
As $\tau^p$ belongs to $L^\infty(\Omega)$, we obtain
\begin{equation}\label{toto2b}
\lim _{\e\to 0} \IOm \tau^p(\vec{x})\beta^p\left(\frac{{\bf R}\vec{x}}{\eta}\right)\,\mathrm{d}\vec{x} = \IOm \tau^p(\vec{x}) \,\mathrm{d}\vec{x} \int_{Y^m} \beta^p(\vec{y}) \mathrm{d}\vec{y}
  \end{equation}
From Fubini's theorem, this implies that (\ref{toto2}) holds.

This result is extended by linearity to step functions $\phi_k\in S_t(\Omega,C_\sharp(Y^m))$ such that
$\phi_k=\displaystyle{\sum_{i=1}^k t_i \chi_{A_i}(\vec{x})\psi_i(\vec{y})}$, where $A_i=\{\vec{x}\in\Omega \, , \, \phi_k(\vec{x},.)=t_i\}$ and $\psi_i(\vec{y})\in C_\sharp(Y^m)$. We deduce that (\ref{toto2}) holds by density in $L^p(\Omega,C_\sharp(Y^m))$.

More precisely, we consider $\phi_k\in S_t(\Omega,C_\sharp(Y^m))$. There exists a sequence of step functions $\phi_k=\displaystyle{\sum_{i=1}^k t_i \chi_{A_i}(\vec{x})\psi_i(\vec{y})}$
such that
\begin{equation}
    \lim_{k\to \infty}\int_\Omega\left(\sup_{\vec{y}\in Y^m}\mid\phi_k\xy-\phi\xy\mid\right)^p \mathrm{d}\vec{x}=\lim_{k\to \infty}{\Vert
\phi_k-\phi\Vert}^p_{L^p(\O,C_\sharp(Y^m))}
\end{equation}
Moreover, from the triangular inequality and the continuity of the linear map ${\bf R}$,
we deduce that there exists a constant $C >0$ such that
\begin{equation}
    {\left\Vert
\phi\xxe\right\Vert}_{L^p(\O)}
\leq C{\left\Vert
\phi\xxe-\phi_k\xxe\right\Vert}_{L^p(\O)}
+{\left\Vert\phi_k\xxe\right\Vert}_{L^p(\O)}
\label{dens1}
\end{equation}

Noting that for every $v\in L^p(\O,C_\sharp(Y^m))$
\begin{equation}
    \int_\Omega {\left| v\xxe \right|}^p \mathrm{d}\vec{x} \leq
    \int_\Omega\left(\sup_{\vec{y}\in Y^m}\left| v\xy\right|\right)^p \mathrm{d}\vec{x}={\Vert
v\Vert}^p_{L^p(\O,C_\sharp(Y^m))}
\end{equation}
we deduce from (\ref{dens1}) that for every integer $k$
\begin{equation}
    {\left\Vert
\phi\xxe\right\Vert}_{L^p(\O)}
\leq C{\Vert
\phi-\phi_k\Vert}_{L^p(\O,C_\sharp(Y^m))}
+{\left\Vert\phi_k\xxe\right\Vert}_{L^p(\O)}
\label{dens2}
\end{equation}
Since $\phi_k$ is admissible, we deduce that there exists a constant $C>0$ such that
\begin{equation}
\begin{array}{ll}
\lims
    {\left\Vert
\phi\xxe\right\Vert}_{L^p(\O)}
&\leq C{\Vert
\phi-\phi_k\Vert}_{L^p(\O,C_\sharp(Y^m))}
+{\Vert\phi_k\Vert}_{L^p(\O\times Y^m)} \\
&\leq C{\Vert
\phi-\phi_k\Vert}_{L^p(\O,C_\sharp(Y^m))}
+{\Vert\phi_k-\phi\Vert}_{L^p(\O\times Y^m)}
+{\Vert\phi\Vert}_{L^p(\O\times Y^m)}\\
&\leq (C+1){\Vert
\phi-\phi_k\Vert}_{L^p(\O,C_\sharp(Y^m))}
+{\Vert\phi\Vert}_{L^p(\O\times Y^m)}
\end{array}
\label{dens4}
\end{equation}
Passing to the limit $k\to  \infty$ in (\ref{dens4}), we obtain
\begin{equation}
\lims
    {\left\Vert
\phi\xxe\right\Vert}_{L^p(\O)}
\leq
{\Vert\phi\Vert}_{L^p(\O\times Y^m)}
\label{dens5}
\end{equation}
Similar arguments hold for the lower limit. In that case, we obtain that for every
integer $k$
\begin{equation}
\begin{array}{ll}
\limf
    {\left\Vert
\phi\xxe\right\Vert}_{L^p(\O)}
&\geq {\Vert\phi_k\Vert}_{L^p(\O\times Y^m)}
-C{\Vert
\phi-\phi_k\Vert}_{L^p(\O,C_\sharp(Y^m))}
\\
&\geq {\Vert\phi\Vert}_{L^p(\O\times Y^m)}
-{\Vert\phi_k-\phi\Vert}_{L^p(\O\times Y^m)}
-C{\Vert
\phi-\phi_k\Vert}_{L^p(\O,C_\sharp(Y^m))}\\
&\geq {\Vert\phi\Vert}_{L^p(\O\times Y^m)}
-(C+1){\Vert
\phi-\phi_k\Vert}_{L^p(\O,C_\sharp(Y^m))}
\end{array}
\label{dens6}
\end{equation}

Passing to the limit $k\to  \infty$ in (\ref{dens6}), we obtain
\begin{equation}
\limf
    {\left\Vert
\phi\xxe \right\Vert}_{L^p(\O)}
\geq
{\Vert\phi\Vert}_{L^p(\O\times Y^m)}
\label{dens7}
\end{equation}
Equality (\ref{toto2}) is established combining inequalities (\ref{dens5}) and (\ref{dens7}).

\hfill \boxed{}

A second Corollary  of proposition \ref{L1convergence} is
\begin{corollary}
\label{L2to2scale}
Let $\O$ be an open bounded set in $\R^n$ and $Y^m={]0,1[}^m$ with $m>n$.
Let a sequence $(u_\e)$ converge strongly  
to $u$ in $L^p(\O)$. Then there exist a vanishing
subsequence $\e_k$ and a limit $\uo \in L^p(\O\times Y^m)$
($Y^m$-periodic in $\vec{y}$) such that $\uek \wdto u_0=u$ as $\e_k \to 0$ for any $m\times n$ matrix ${\bf R}$ satisfying \eqref{criterion}.   
\end{corollary}

{\bf Proof}  
Let  $\phi(\vec{x},\vec{y})\in L^q(\O; C_\sharp(Y^m))$. Then, we have
\begin{equation}
\begin{array}{ll}
&\left|\displaystyle{\IOm\ue (\vec{x}) \phi\left(\vec{x}, \frac{{\bf R}\vec{x}}{\e}\right) \,\mathrm{d}\vec{x}
-
\IOAY \uo \phi(\vec{x},\vec{y})\,\mathrm{d}\vec{x}\mathrm{d}\vec{y}} \right| \\
&\leq
{\Vert u_\e(\vec{x})-u(\vec{x})\Vert}_{L^p(\O)}
{\left\Vert \phi\left(\vec{x}, \frac{{\bf R}\vec{x}}{\e}\right)\right\Vert}_{L^q(\O)} \\
&+\left|\displaystyle{\IOm u(\vec{x}) \phi\left(\vec{x}, \frac{{\bf R}\vec{x}}{\e}\right) \,\mathrm{d}\vec{x}
-
\IOAY u(\vec{x}) \phi(\vec{x},\vec{y})\,\mathrm{d}\vec{x}\mathrm{d}\vec{y}} \right|
\end{array}
\end{equation}
Since $u(\vec{x}) \phi(\vec{x},\vec{y})\in L^1(\O; C_\sharp(Y^m))$ as $u$ is in $L^p(\Omega)$ by assumption, proposition~\ref{L1convergence} ensures
the second term in RHS goes to zero. The first term in RHS goes to zero since
by corollary~\ref{Lpconvergence}
\begin{equation}\label{toto2c}
\lim _{\e\to 0} \IOm {\left| \phi\left(\vec{x}, \frac{{\bf R}\vec{x}}{\e}\right)\right|}^q \,\mathrm{d}\vec{x} = \IOAY {\mid\phi(\vec{x},\vec{y})\mid}^q\,\mathrm{d}\vec{x} \mathrm{d}\vec{y}
  \end{equation}
and $(u_\e)$ converges strongly to $u$ in $L^p(\O)$ by assumption.
\hfill \boxed{}

The following proposition provides us with
a corrector type result for the sequence
$u_\e$ when its limit $u_0$ is smooth enough:

\begin{proposition}\label{allaireLp}
Let ${\bf R}$ be a linear map from $\R^n$ in $\R^m$ satisfying
(\ref{criterion}). Let $\ue$ be a sequence
such that $\ue\wdto \uo$ (weakly). Then

i) $\ue$ weakly converges in $L^p(\O)$ towards $\displaystyle{u(\vec{x}) = \int_{Y^m} u_0(\vec{x},\vec{y}) \;\mathrm{d}\vec{y}}$ and
\begin{equation}
\limf {\Vert\ue\Vert}_{L^p(\O)}\geq {\Vert u_0 \Vert}_{L^p(\O \times
Y^m)} \geq {\Vert u\Vert}_{L^p(\O)}  
\end{equation}
ii) Let $\ve$ be another bounded sequence in $L^q(\O)$, $1/p + 1/q = 1$, such that
$v_\e\dto v_0$ (strongly), then
\begin{equation}
\ue\ve\to  w(\vec{x}) \; \hbox{ in ${\mathcal{D}}'(\O)$ where
$\displaystyle{w(\vec{x}) = \int_{Y^m}\uo v_0(\vec{x},\vec{y})
\,\mathrm{d}\vec{y}}$}
\end{equation} 
iii) If ${\displaystyle\ue\dto u_0(\vec{x},\vec{y})}$ strongly and $u_0$ is smooth enough (\eg $u_0 \in L^p(\Omega, C_{\sharp}(Y^m))$) and  
\begin{equation}
\displaystyle {{\left\Vert u_0\xxe \right\Vert}_{L^p(\O)} \to{\Vert
u_0\xy\Vert}_{L^p(\O\times Y^m)}}
\end{equation} then
\begin{equation}
\displaystyle{{\left\| \ue - u_0\xxe \right\|}_{L^p(\O)}\to  0}
\label{3i}
\end{equation}
\end{proposition}

{\bf Proof } 
(i) Choosing test functions $\phi$ in $L^q(\O; C_\sharp(Y^m))$ independent of the $y$ variable in definition \ref{2scaleweak}, one has that for every $\phi\in L^q(\O)$
  \begin{equation}\label{eq:weak2scaleagain}
  \begin{array}{ll}
\displaystyle{\lim _{\e\to 0} \IOm\ue (\vec{x}) \phi\left(\vec{x}\right) \,\mathrm{d}\vec{x}}
&= \displaystyle{\IOAY \uo \phi(\vec{x})\,\mathrm{d}\vec{x}\mathrm{d}\vec{y}} \\
&= \displaystyle{\int_\Omega \phi\left(\vec{x}\right) \left( \int_{Y^m} \uo \mathrm{d}\vec{y} \right) \, \mathrm{d}\vec{x}}
\end{array}
  \end{equation}
Moreover, $(u_\e)$ is bounded in $L^p(\O)$ as a weakly convergent sequence in a Hilbert space.

Then, let $\varphi_m$ be a sequence in $L^q(\O,C_\sharp(Y^m))$ such that
$\varphi_m$ converges to $\mid u_0 \mid^{p-2} u_0$ strongly in $L^q(\O\times Y^m)$.

We first apply the Young inequality for real numbers $a$ and $b$, and $1<p<\infty$, $1/p+1/q=1$, which states that
$ab\leq \mid a\mid^p/p + \mid b \mid^q/q$. We consider $a=u_\eta$ and $b=\varphi_m$
to get
\begin{equation}
    \IOm \mid {u}_{\e}(\vec{x})\mid^p \,\mathrm{d}\vec{x}\geq p\IOm {u}_{\e}(\vec{x}) {\varphi}_m\left(\vec{x},\frac{ {\bf
R}\vec{x}}{\e}\right) \,\mathrm{d}\vec{x}-(p-1)\IOm \left| {\varphi}_m\left(\vec{x},\frac{ {\bf
R}\vec{x}}{\e}\right)\right|^q \,\mathrm{d}\vec{x}
\end{equation}
We first pass to the limit when $\eta$ goes to zero:
\begin{equation}
\limf {\Vert\ue\Vert}^p_{L^p(\O)}\geq
p\IOAY   {u}_0(\vec{x},\vec{y}) {\varphi}_m(\vec{x},\vec{y})\,\mathrm{d}\vec{x}\mathrm{d}\vec{y}
-(p-1)\IOAY {\left|{\varphi}_m(\vec{x},\vec{y})\right|}^q\,\mathrm{d}\vec{x}\mathrm{d}\vec{y}
\end{equation}
where we have used that $\ue\wdto \uo$ weakly and ${\displaystyle{\varphi}_m\left(\vec{x},\frac{ {\bf
R}\vec{x}}{\e}\right)\dto {\varphi}_m(\vec{x},\vec{y})}$ strongly (making use of proposition \ref{Lpconvergence}). 

We then pass to the limit when $m$ goes to infinity:
\begin{equation}
\limf {\Vert\ue\Vert}^p_{L^p(\O)}\geq
p\IOAY {\mid{u}_0(\vec{x},\vec{y})\mid}^p\,\mathrm{d}\vec{x}\mathrm{d}\vec{y}
-(p-1)\IOAY {\mid{u}_0(\vec{x},\vec{y})\mid}^p\,\mathrm{d}\vec{x}\mathrm{d}\vec{y}={\Vert u_0 \Vert}^p_{L^p(\O \times
Y^m)}
\label{young}
\end{equation}
where we have used that $\varphi_m\to\mid u_0 \mid^{p-2} u_0$
strongly in $L^q(\O\times Y^m)$.

Moreover, thanks to Jensen's inequality, we have that
\begin{equation}
{\Vert u \Vert}^p_{L^p(\O)}
=
\int_\Omega \left| \int_{Y^m} {{u}_0(\vec{x},\vec{y})}\mathrm{d}\vec{y} \right|^p\,\mathrm{d}\vec{x}
\leq
\IOAY {\mid{u}_0(\vec{x},\vec{y})\mid}^p\,\mathrm{d}\vec{x}\mathrm{d}\vec{y}={\Vert u_0 \Vert}^p_{L^p(\O\times Y^m)}
\label{jensen}
\end{equation}
We conclude by combining (\ref{young}) and (\ref{jensen}).

(ii) Let $\psi_m$ be a sequence in $L^p(\O,C_\sharp(Y^m))$ such that
$\psi_m$ converges to $u_0$ strongly in $L^p(\O\times Y^m)$.   Let $\tau$ be a function in $C^{\infty}_0(\Omega)$.   We note that $\ve\tau\dto v_0\tau$ (strongly). Thus, passing to the two-scale-cut-and-projection limit when $\e$ goes to zero in the product of $v_\e\tau$ and $\psi_m$, we have that
\begin{equation}
    \lim_{\eta\to 0}\IOm {\psi}_m\left(\vec{x},\frac{ {\bf
R}\vec{x}}{\e}\right)v_\e(\vec{x})\tau(\vec{x}) \,\mathrm{d}\vec{x}
=\IOAY {\psi}_m\left(\vec{x},\vec{y}\right)v_0(\vec{x},\vec{y})\tau(\vec{x})\,\mathrm{d}\vec{x}\mathrm{d}\vec{y}
\end{equation}
We then pass to the limit when $m$ goes to infinity
\begin{equation}
    \lim_{m\to\infty}\lim_{\eta\to 0}\IOm {\psi}_m\left(\vec{x},\frac{ {\bf
R}\vec{x}}{\e}\right)v_\e(\vec{x})\tau(\vec{x}) \,\mathrm{d}\vec{x}
=\IOAY u_0\left(\vec{x},\vec{y}\right)v_0(\vec{x},\vec{y})\tau(\vec{x})\,\mathrm{d}\vec{x}\mathrm{d}\vec{y}
\label{trick1}
\end{equation}
where we have used that $\psi_m$ converges to $u_0$ strongly in $L^p(\O\times Y^m)$.

Moreover, from the triangular inequality, we have
\begin{equation}
\begin{array}{ll}
&\displaystyle{\left|\IOm u_\e(\vec{x})v_\e(\vec{x})\tau(\vec{x}) \,\mathrm{d}\vec{x}
-\IOAY u_0\left(\vec{x},\vec{y}\right)v_0(\vec{x},\vec{y})\tau(\vec{x})\,\mathrm{d}\vec{x}\mathrm{d}\vec{y}\right|} \\
&\leq \displaystyle{\left|\IOm \left[u_\e(\vec{x})-{\psi}_m\left(\vec{x},\frac{ {\bf
R}\vec{x}}{\e}\right)\right]v_\e(\vec{x})\tau(\vec{x}) \,\mathrm{d}\vec{x}\right|} \\
&+\displaystyle{\left|\IOm {\psi}_m\left(\vec{x},\frac{ {\bf
R}\vec{x}}{\e}\right)v_\e(\vec{x})\tau(\vec{x}) \,\mathrm{d}\vec{x}
-\IOAY u_0\left(\vec{x},\vec{y}\right)v_0(\vec{x},\vec{y})\tau(\vec{x})\,\mathrm{d}\vec{x}\mathrm{d}\vec{y}\right|}
\end{array}
\label{trick2}
\end{equation}
Combining (\ref{trick1}) and (\ref{trick2}) we get
\begin{equation}
\begin{array}{ll}
&\displaystyle{\lims\left|\IOm u_\e(\vec{x})v_\e(\vec{x})\tau(\vec{x}) \,\mathrm{d}\vec{x}
-\IOAY u_0\left(\vec{x},\vec{y}\right)v_0(\vec{x},\vec{y})\tau(\vec{x})\,\mathrm{d}\vec{x}\mathrm{d}\vec{y}\right|} \\
&\leq \displaystyle{\limsup_{m\to\infty}\lims \left|\IOm \left[u_\e(\vec{x})-{\psi}_m\left(\vec{x},\frac{ {\bf
R}\vec{x}}{\e}\right)\right]v_\e(\vec{x})\tau(\vec{x}) \,\mathrm{d}\vec{x} \right|}
\end{array}
\label{trick3}
\end{equation}
It remains to prove that the right-hand side in (\ref{trick3}) vanishes.

We first recall that $\tau\in C^\infty_0(\Omega)$ and we then invoke Hölder's inequality
to deduce that
\begin{equation}
\begin{array}{ll}
&\displaystyle{\left|\IOm \left[u_\e(\vec{x})-{\psi}_m\left(\vec{x},\frac{ {\bf
R}\vec{x}}{\e}\right)\right]v_\e(\vec{x})\tau(\vec{x}) \,\mathrm{d}\vec{x}\right|} \\
&\leq\max_{\vec{x}\in\Omega}\mid\tau(\vec{x})\mid \, \displaystyle{\left|\IOm \left[u_\e(\vec{x})-{\psi}_m\left(\vec{x},\frac{ {\bf
R}\vec{x}}{\e}\right)\right]v_\e(\vec{x}) \,\mathrm{d}\vec{x}\right|} \\
&\leq\max_{\vec{x}\in\Omega}\mid\tau(\vec{x})\mid \, \displaystyle{\left(\IOm {\left| u_\e(\vec{x})-{\psi}_m\left(\vec{x},\frac{ {\bf
R}\vec{x}}{\e}\right)\right|}^p \,\mathrm{d}\vec{x}\right)^{1/p}}
\displaystyle{\left(\IOm {\mid v_\e(\vec{x})\mid}^q \,\mathrm{d}\vec{x}\right)^{1/q}} \\
&\leq C(\Omega) \, \displaystyle{{\left\| \ue - \psi_m\xxe \right\|}_{L^p(\O)}} 
\end{array}
\label{trick4}
\end{equation}
where in the last inequality we have used that  $v_\e$ is a bounded sequence in $L^q(\Omega)$.

We now invoke the Clarkson inequalities applied to functions $u_\e$ and $\psi_m$ in $L^p(\Omega)$:
\begin{equation}
\begin{array}{ll}
{\displaystyle \frac{1}{2^p}\left\|{u_\e-\psi_m\xxe}\right\|_{L^{p}(\Omega)}^{p}}&\displaystyle{ \leq {\frac {1}{2}}\left(\|u_\e\|_{L^{p}(\Omega)}^{p}+\|\psi_m\xxe\|_{L^{p}(\Omega)}^{p}\right)} \\
&\displaystyle{-\left\|{\frac {u_\e+\psi_m\xxe}{2}}\right\|_{L^{p}(\Omega)}^{p} \mbox{ for }1 < p < 2}
\end{array}
\label{clarksson3}
\end{equation}
and
\begin{equation}
\begin{array}{ll}
     {\displaystyle \frac{1}{2^q}\left\|{u_\e-\psi_m\xxe}\right\|_{L^{p}(\Omega)}^{q}}&\leq \displaystyle{\left({\frac {1}{2}}\|u_\e\|_{L^{p}(\Omega)}^{p}+{\frac {1}{2}}\|\psi_m\xxe\|_{L^{p}(\Omega)}^{p}\right)^{\frac {q}{p}}} \\
     &\displaystyle{-\left\|{\frac {u_\e+\psi_m\xxe}{2}}\right\|_{L^{p}(\Omega)}^{q} \mbox{ for } p\geq 2}
\end{array}     
     \label{clarksson4}
\end{equation}

Passing to the 2-scale cut-and-projection limit in (\ref{clarksson3})
\begin{equation}
\begin{array}{ll}
\lims{\displaystyle \left\|{u_\e-\psi_m\xxe}\right\|_{L^{p}(\Omega)}^{p}}
&\leq \displaystyle{2^{p-1}\left(\|u_0\xy\|_{L^{p}(\Omega\times Y^m)}^{p}+\|\psi_m\xy\|_{L^{p}(\Omega\times Y^m)}^{p}\right)} \\
&-\displaystyle{2^p\left\|{\frac {u_0\xy+\psi_m\xy}{2}}\right\|_{L^{p}(\Omega\times Y^m)}^{p} \mbox{ for }1 < p < 2}
\end{array}
\label{clarksson5}
\end{equation}
and similarly for (\ref{clarksson4})

\begin{equation}
\begin{array}{ll}
\lims{\displaystyle \left\|{u_\e-\psi_m\xxe}\right\|_{L^{p}(\Omega)}^{q}}
&\leq \displaystyle{2^{}\left(\|u_0\xy\|_{L^{p}(\Omega\times Y^m)}^{p}+\|\psi_m\xy\|_{L^{p}(\Omega\times Y^m)}^{p}\right) ^{\frac{1}{p-1}}}\\
&-\displaystyle{ 2^{q}\left\|{\frac {u_0\xy+\psi_m\xy}{2}}\right\|_{L^{p}(\Omega\times Y^m)}^{q} \mbox{ for }  p \geq 2}
\end{array}
\label{clarksson6}
\end{equation}
where we have used
also that $1/p+1/q=1$ thus $q/p=1/(p-1)$.

We now note that
$\displaystyle \left\|{\frac {u_0\xy+\psi_m\xy}{2}}\right\|_{L^{p}(\Omega\times Y^m)}\to\displaystyle \left\|{u_0\xy}\right\|_{L^{p}(\Omega\times Y^m)}$ when $m$ tends to infinity. Thus, taking the lim sup on $m$ in both sides of (\ref{clarksson5}) and (\ref{clarksson6}), we are ensured that
\begin{equation}
\limsup_{m\to\infty}\lims{\displaystyle \left\|{u_\e-\psi_m\xxe}\right\|_{L^{p}(\Omega)}}=0 \mbox{ for } 1<p<\infty
\label{clarksson7}
\end{equation}
We have thus proved that the RHS of (\ref{trick3}) vanishes combining (\ref{trick4}) and (\ref{clarksson7}).

(iii) We finally want to show (\ref{3i}). Assuming that $u_0 \in L^p(\Omega, C_{\sharp}(Y^m))$, we can replace $\psi_m$ by $u_0$ in the Clarkson inequalities (\ref{clarksson3}) and (\ref{clarksson4}), which leads to \begin{equation}
\limsup_{m\to\infty}\lims{\displaystyle \left\|{u_\e-u_0\xxe}\right\|_{L^{p}(\Omega)}}=0 \mbox{ for } 1<p<\infty
\label{clarksson7bis}
\end{equation}
which is (\ref{3i}).

\hfill \boxed{}

\begin{remark}
Our proof of proposition \ref{allaireLp} follows closely that of
Theorems 10 and 11 in \cite{lukkassen2002two}.

We stress that some regularity is needed for $u_0$ in (iii) of proposition \ref{allaireLp}. We refer to Theorem 12 in \cite{lukkassen2002two} for a related result on regularity of two-scale limit that remains valid for two-scale-cut-and-projection limit (its proof is a straightforward extension of that in \cite{lukkassen2002two}, but is rather technical and lengthy).  
\end{remark}

\begin{remark}
We point out that Proposition~\ref{allaireLp} i) does not hold if weak two-scale convergence $\uek \wdto u_0$ is replaced by distributional two-scale convergence (see definition \ref{2scaledistribution}).
Indeed, the choice of space $L^q(\O,C_\sharp(Y^m))$ for test functions $\varphi$ in (\ref{2scaledistribution}) is essential.
This is exemplified by the counter-example in \cite{lukkassen2002two} of a sequence $(u_\eta)$ in $L^p(0,1)$ defined by
$u_\eta(x)=1/\eta$ if $0<x<\eta$ and $u_\eta(x)=0$ if $\eta<x<1$.
One can see that \eqref{eq:distr2scale} is satisfied for $\varphi\in{\mathcal{D}}((0,1);C^\infty_\sharp(Y^m))$ and a two-scale limit $u_0(x,y)=0$. However, considering the test function $g(x)=1$ which is in $L^q(0,1)$, we get $\lim_{\eta\to 0} \int_0^1 u_\e(x)g(x) dx=1$, and so $(u_\eta)$ does not converge to $u_0(x,y)=0$ weakly in $L^p(0,1)$.
\end{remark}

Classes of functions such that $\displaystyle {{\left\| u_0 \xxe
\right\|}_{L^p(\O)} \to{\left\Vert u_0\xy\right\Vert}_{L^p(\O\times Y^m)}}$ are
said to be admissible for the two-scale (cut-and-projection)
convergence. In particular, classes of functions in
$L^p(\O,C_\sharp(Y^m))$ (dense subset in $L^p(\O\times Y^m)$) are
admissible.

In order to homogenize nonlinear PDEs with a monotone partial differential operator as in (\ref{eq:NL_static}), we need to identify the differential relationship between
$\vec{\chi}$ and $u_0$, given a bounded sequence $(u_\e)$ in $W^{1,p}(\O)$
(such that $u_\e\wdto u_0$ and $\nabla u_\e\wdto \vec{\chi}$).
This problem was solved by
Allaire in the case of periodic functions \cite{Allaire1992}
and extended by Bouchitt\'e et al. for quasiperiodic functions
\cite{Bouchitte+etal2010} in $W^{1,2}(\O)$ and  revisited in \cite{Wellander+etal2017,Wellander+etal2019}.

\section{Function spaces for cut-and-projection partial differential operators}  
To carry out the homogenization analysis of nonlinear PDEs defined on quasiperiodic domains, we need to pass to the limit  when $\eta$ goes to zero in gradient and 
divergence operators acting on solutions of PDEs.  To do this we introduce some suitable function spaces and for this we will define differential operators acting on the $\R^n$-subspace  in  $\R^m$.  

Defined as in \cite{Wellander+etal2017} they are given as
\begin{equation*}
 \nabla_{\bf R}\;u(\vec{y}) =
 \mathrm{grad}_{\bf R}\;u(\vec{y}) = {\bf R}^T \mathrm{grad}_{\vec{y}} \;  u(\vec{y}) = {\bf R}^T \nabla_{\vec{y}} \;  u(\vec{y}) 
\end{equation*}
\begin{equation*}
\mathrm{div}_{\bf R}\;\vec{u}(\vec{y}) = 
{\bf R}^T \nabla_{\vec{y}} \cdot \vec{u}(\vec{y}) 
\end{equation*} 
We define the following functions spaces associated with the differential operators defined above
\begin{equation}
\begin{array}{ll}
\displaystyle
{\mathcal{W}}^p_{\sharp}(\mathrm{grad}_{\bf R},Y^m) = & \displaystyle{\Bigl \{u\in
{L^p_\sharp(Y^m)} \; \mid \; \mathrm{grad}_{\bf R}\;u \in
{L^p_\sharp(Y^m; \R^n)}
\Bigr \}
}
\end{array}
\label{def:grad-Yspace}
\end{equation}
\begin{equation}
\begin{array}{ll}
\displaystyle
{\mathcal{W}}^p_{\sharp}(\mathrm{div}_{\bf R}, Y^m) = & \displaystyle{\Bigl \{ \vec{u}\in
{L^p_\sharp(Y^m; \R^n)} \; \mid \; \mathrm{div}_{\bf R}\;\vec{u}  \in
{L^p_\sharp(Y^m)}
\Bigr \}
}
\end{array}
\label{def:divR-Yspace}
\end{equation}

and 
 \begin{equation}
  {\mathcal{W}^p}(\mathrm{div},\Omega) =  \displaystyle{\Bigl \{ \vec{u}\in
{L^p(\Omega; \R^n)} \;\mid \; \mathrm{div} \;\vec{u} \in {L^p(\Omega}   \Bigr \}
}
 \end{equation}

 We have the following integration by parts type generalization to the $L^p$ case  of Lemma~6 given in the $L^2$ setting in \cite{Wellander+etal2019}

\begin{lemma}[Green's Identity] \label{lem:Integration_by_parts_green}
It holds that
\begin{equation}\label{eq:Green_identity}
 -\int_{Y^m}  \left({\bf R} {\bf R}^T \nabla_y \right) \cdot \vec{\phi}(\vec{y}) \;   {\theta}(\vec{y}) \;\mathrm{d}\vec{y} = \int_{Y^m}    \vec{\phi}(\vec{y})  \cdot  \left({\bf R} {\bf R}^T \nabla_y \right)   \;{\theta}(\vec{y}) \;\mathrm{d}\vec{y}
\end{equation}
for  $\vec{\phi} \in  
{\mathcal{W}}^q_{\sharp}(\mathrm{div}_{\bf R}, Y^m)$ and  $\theta \in
 {\mathcal{W}}^p_{\sharp}(\mathrm{grad}_{\bf R},Y^m)$, $1/p + 1/q =1$.
\end{lemma}
{\bf Proof }
The proof relies on standard matrix operations and the well known extension from $W^{1,2}$ Sobolev spaces to the setting of $W^{1,p}$ and $W^{1,q}$-duality pairing 
\cite{Brezis2010}. The periodic boundary conditions imply
\begin{equation}
\begin{aligned}
&
 -\int_{Y^m}  \left({\bf R} {\bf R}^T \nabla_y \right) \cdot \vec{\phi}(\vec{y})  \;   {\theta}(\vec{y}) \;\mathrm{d}\vec{y} =
 -\int_{Y^m}  \left({\bf R}^T \nabla_y \right) \cdot {\bf R}^T \vec{\phi}(\vec{y})  \;   {\theta}(\vec{y}) \;\mathrm{d}\vec{y} = \\
 &
 -\int_{Y^m}  \nabla_y \cdot {\bf R}{\bf R}^T \vec{\phi}(\vec{y})  \;   {\theta}(\vec{y}) \;\mathrm{d}\vec{y} =  
 \int_{Y^m}   {\bf R} {\bf R}^T \vec{\phi} \cdot  \nabla_y \;{\theta}(\vec{y}) \;\mathrm{d}\vec{y}=
 \int_{Y^m}    \vec{\phi}(\vec{y})  \cdot {\bf R} {\bf R}^T\nabla_y \;{\theta}(\vec{y}) \;\mathrm{d}\vec{y}
\end{aligned}
\end{equation}
for any pair of functions for  $\vec{\phi} \in  
{\mathcal{W}}^q_{\sharp}(\mathrm{div}_{\bf R}, Y^m)$ and  $\theta \in
 {\mathcal{W}}^p_{\sharp}(\mathrm{grad}_{\bf R},Y^m)$, $1/p + 1/q =1$.
\hfill \boxed{}

This motivates us to make the following decomposition as in the $L^2$-case in \cite{Wellander+etal2019}. We  decompose 
$W^{1,p}_\sharp(Y^m)$ into two spaces,
\begin{equation}\label{splitXp}
W^{1,p}_\sharp(Y^m) = X_p \oplus  X_p^\perp
\end{equation}
where
\begin{equation}\label{eq:W12_composition}
    X_p^\perp = \left\{ u \in W^{1,p}_\sharp(Y^m) | {\bf R}{\bf R}^T \; \nabla_y u = \vec{0} \right\}
\end{equation}
and
\begin{equation}\label{eq:X-proj_space}
X_p = \left\{ u \in W^{1,p}_\sharp(Y^m) | \left({\bf I}_m - {\bf R}{\bf R}^T\right) \; \nabla_y u = \vec{0} \right\}
\end{equation} 


\begin{remark}
Note that the projection of a vector $\vec{v}$ in $\R^m$ on $\R^n$, ${\bf R}^T \; \vec{v} = \vec{0}$, where $\vec{0}$ is the zero vector in $\R^n$, implies that $\vec{v}$ is orthogonal to the hyperplane, \ie orthogonal to $\R^n$.  It follows that ${\bf R}{\bf R}^T \; \vec{v} = \vec{0}$, where $\vec{0}$ is the zero vector in $\R^m$. Vectors $\vec{w}$ in $\R^m$, orthogonal to $\vec{v}$, satisfies $\left({\bf I}_m - {\bf R}{\bf R}^T \right)\; \vec{w} = \vec{0}$.  We conclude that $X_p$ contains all functions in $W^{1,p}_\sharp(Y^m)$ whose gradients have  all their components in the hyperplane, $\R^n$, which means that  $X_p$ can be identified with  ${\mathcal{W}}^p_{\sharp}(\mathrm{grad}_{\bf R},Y^m) $.
\end{remark}

\section{Compactness, homogenization and corrector results}

\begin{proposition}\label{prop:grad-split}
Let  $\{u_\e\}$ be a uniformly bounded
sequence in $W^{1,p}(\Omega)$.
Then there exist a subsequence $\{u_{\e_k}\}$  
 and  functions $u\in W^{1,p}(\Omega)$ and $ {u}_1\xy \in L^p(\Omega, X_p )$ 
such that
\begin{equation}\label{eq:gradsplit_2s_strong}
u_{\e_k} \dto u(\vec{x}),  \qquad \mathrm{grad}\; u_{\e_k} \wdto \mathrm{grad} \; u (\vec{x}) + \mathrm{grad}_{\bf R} \;  u_1\xy ,  \qquad 
{\e_k} \to 0 
\end{equation}
\end{proposition}  
\begin{remark}
Note that ${u}_1\xy \in L^p(\Omega, X_p)$ implies that   $\mathrm{grad}_{\bf R} \;  {u}_1\xy \in L^p(\Omega, L^p_\sharp(Y^m; \R^n ))$ and that we don't have enough regularity  to ensure  that  $ {u}_1$  belongs to $ L^p(\Omega, W^{1,p}_\sharp(Y))$. Indeed, we cannot say anything about the regularity of  $u_1$ in the  direction orthogonal to the hyperplane.  However, in the decomposition in  (\ref{splitXp}), the gradient of the potential in the "direction " of the hyperplane can be obtained from the gradient of the potential via a rotation of coordinate system.
Further note that unlike in \cite{Wellander+etal2017,Wellander+etal2019}, in the proof below we use the notion of 2-scale cut-and-projection convergence in distributional sense (Definition \ref{2scaledistribution}), and not in weak sense (Definition \ref{2scaleweak}).  
\end{remark}

{\bf Proof }
The first assertion follows by the compact embedding of $L^p(\Omega)$ in $ W^{1,p}(\Omega)$,  Propositions~\ref{prop:weaktwoscaleprop} and Definition~\ref{def:strong_two_scale}. Note that $ {\bf R} \nabla {u}_{\e}(\vec{x}) $, $\vec{x} \in \Omega$ is uniformly  bounded in $L^p(\Omega;\R^m)$.
Let $\vec{\varphi}\in{\mathcal{D}}(\O;C^\infty_\sharp(Y^m))^m$.
We have the following identities 
 \begin{equation}\label{eq:grad-split-compactness_proof_0}
  \begin{aligned}
 &
 \lim _{\e\to 0} \IOm    {\bf R} \nabla {u}_{\e}(\vec{x}) \cdot \vec{\varphi}\left(\vec{x},\frac{ {\bf
R}\vec{x}}{\e}\right) \,\mathrm{d}\vec{x} = 
\lim _{\e\to 0} \IOm     \nabla  {u}_{\e}(\vec{x}) 
\cdot {\bf R}^T \vec{\varphi}\left(\vec{x},\frac{ {\bf
R}\vec{x}}{\e}\right) \,\mathrm{d}\vec{x}= \\
&
-\lim _{\e\to 0} \IOm     {u}_{\e}(\vec{x}) 
 \left(\nabla_x \cdot {\bf R}^T \vec{\varphi}\left(\vec{x},\frac{ {\bf
R}\vec{x}}{\e}\right) +
\e^{-1}\left({\bf R}^T \nabla_y\right) \cdot {\bf R}^T \vec{\varphi}\left(\vec{x},\frac{ {\bf
R}\vec{x}}{\e}\right)\right) \,\mathrm{d}\vec{x}  
 \end{aligned}
 \end{equation}
and
 \begin{equation}\label{eq:local_divergence_of_testf}
   \left( {\bf R}^T \nabla_y \right) \cdot
    {\bf R}^T \vec{\varphi}\left(\vec{x},\frac{ {\bf R}\vec{x}}{\e}\right)
    =  \left({\bf R} {\bf R}^T \nabla_y \right) \cdot \vec{\varphi}\left(\vec{x},\frac{ {\bf R}\vec{x}}{\e}\right)
 \end{equation}
 Multiplying both sides in \eqref{eq:grad-split-compactness_proof_0}  with $\e$ and Lemma~\ref{lem:Integration_by_parts_green} gives  the limit
 \begin{equation}\label{eq:local_divergence_of_testf_1}
  \begin{aligned}
  0 = &  \IOAY  u_0(\vec{x},\vec{y}) \left( {\bf R}^T \nabla_y \right) \cdot {\bf R}^T \vec{\varphi}\left(\vec{x},\vec{y}\right)
  \,\mathrm{d}\vec{x}\mathrm{d}\vec{y}
  =  \IOAY  u_0(\vec{x},\vec{y}) \left( {\bf R} {\bf R}^T \nabla_y \right) \cdot  \vec{\varphi}\left(\vec{x},\vec{y}\right)
  \,\mathrm{d}\vec{x}\mathrm{d}\vec{y}=\\
  &
  \IOAY  {\bf R} {\bf R}^T \nabla_y u_0(\vec{x},\vec{y}) \cdot  \vec{\varphi}\left(\vec{x},\vec{y}\right)\,\mathrm{d}\vec{x}\mathrm{d}\vec{y}
 \end{aligned}
 \end{equation}
 for all $\vec{\varphi}\in{\mathcal{D}}(\O;C^\infty_\sharp(Y^m))^m$. The interpretation of \eqref{eq:local_divergence_of_testf_1} is that $u_0 \in X_p^\perp$, \ie   the gradient $\nabla_y u_0(\vec{x},\vec{y})$ has no component in the hyper plane in $\R^m$ defined by  ${\bf R}: \R^n \to \R^m$.  Indeed, we  conclude that the potential $ u_0(\vec{x}, \vec{y}) = u(\vec{x}) $, is a function of $\vec{x}$ only  due to the compact embedding of  $L^p(\Omega)$ in $ W^{1,p}(\Omega)$ and that the two-scale limit equals the strong limit, if it exists.
 Next, let 
$\vec{\varphi}\in{\mathcal{D}}(\O;C^\infty_\sharp(Y^m))^m$, and $\vec{\psi}\in {\mathcal{D}}(\O;C^\infty_\sharp(Y^m))^n$.  We have three  limits 
   \begin{equation}\label{eq:grad-split-compactness_proof_2}
  \begin{aligned}
 &
 \lim _{\e\to 0} \IOm    {\bf R} \nabla {u}_{\e}(\vec{x}) \cdot \vec{\varphi}\left(\vec{x},\frac{ {\bf
R}\vec{x}}{\e}\right) \,\mathrm{d}\vec{x} =  
\IOAY   \hat{\vec{\chi}}_0(\vec{x},\vec{y}) \cdot  \vec{\varphi}(\vec{x},\vec{y})\,\mathrm{d}\vec{x}\mathrm{d}\vec{y}   
 \end{aligned}
 \end{equation}
  \begin{equation}\label{eq:grad-split-compactness_proof_1}
  \begin{aligned}
 &
 \lim _{\e\to 0} \IOm    {\bf R} \nabla {u}_{\e}(\vec{x}) \cdot \vec{\varphi}\left(\vec{x},\frac{ {\bf
R}\vec{x}}{\e}\right) \,\mathrm{d}\vec{x} = 
\lim _{\e\to 0} \IOm     \nabla  {u}_{\e}(\vec{x}) 
\cdot {\bf R}^T \vec{\varphi}\left(\vec{x},\frac{ {\bf
R}\vec{x}}{\e}\right) \,\mathrm{d}\vec{x}= \\
&
\IOAY \vec{\chi}_0(\vec{x},\vec{y}) \cdot {\bf R}^T  \vec{\varphi}(\vec{x},\vec{y})\,\mathrm{d}\vec{x}\mathrm{d}\vec{y}   =   
\IOAY  {\bf R}\chi_0(\vec{x},\vec{y}) \cdot  \vec{\varphi}(\vec{x},\vec{y})\,\mathrm{d}\vec{x}\mathrm{d}\vec{y}  
 \end{aligned}
 \end{equation}
 and
\begin{equation}\label{eq:grad-split-compactness_proof}
  \begin{aligned}
&
\lim _{\e\to 0} \IOm     \nabla {u}_{\e}(\vec{x}) \cdot \vec{\psi}\left(\vec{x},\frac{ {\bf
R}\vec{x}}{\e}\right) \,\mathrm{d}\vec{x} = 
\lim _{\e\to 0} \IOm     \nabla  {u}_{\e}(\vec{x}) 
\cdot  \vec{\psi}\left(\vec{x},\frac{ {\bf
R}\vec{x}}{\e}\right) \,\mathrm{d}\vec{x}= \\
&
\IOAY \tilde{\vec{\chi}}_0(\vec{x},\vec{y}) \cdot   \vec{\psi}(\vec{x},\vec{y})\,\mathrm{d}\vec{x}\mathrm{d}\vec{y}   \\     
 \end{aligned}
 \end{equation}
 We find that $\hat{\vec{\chi}}_0(\vec{x},\vec{y})=  {\bf R}\vec{\chi}_0(\vec{x},\vec{y}) $ and $\tilde{\vec{\chi}}_0(\vec{x},\vec{y})=  {\bf R}^T\hat{\vec{\chi}}_0(\vec{x},\vec{y}) = \vec{\chi}_0 (\vec{x},\vec{y}) $. Next, choosing   test functions    $\vec{\varphi}\in{\mathcal{D}}(\O;C^\infty_\sharp(Y^m))^m$, such that $\left({\bf R} {\bf R}^T \nabla_y \right)\cdot \vec{\varphi}(\vec{x},\vec{y}) = 0$ in \eqref{eq:grad-split-compactness_proof_0} with 
  \eqref{eq:local_divergence_of_testf}  gives 
  \begin{equation}\label{eq:grad-split-compactness_proof_0_2}
  \begin{aligned}
 &
 \lim _{\e\to 0} \IOm    {\bf R} \nabla {u}_{\e}(\vec{x}) \cdot \vec{\varphi}\left(\vec{x},\frac{ {\bf
R}\vec{x}}{\e}\right) \,\mathrm{d}\vec{x} =   
-\lim _{\e\to 0} \IOm     {u}_{\e}(\vec{x}) 
 \left(\nabla_x \cdot {\bf R}^T \vec{\varphi}\left(\vec{x},\frac{ {\bf
R}\vec{x}}{\e}\right)  \right) \,\mathrm{d}\vec{x}  =  \\
&
- \IOAY u(\vec{x}) 
 \left(\nabla_x \cdot {\bf R}^T \vec{\varphi}\left(\vec{x},\vec{y}\right)  \right) \,\mathrm{d}\vec{x}\mathrm{d}\vec{y} = 
\IOAY  {\bf R} \nabla  u(\vec{x}) 
  \cdot\vec{\varphi}\left(\vec{x},\vec{y}\right)   \,\mathrm{d}\vec{x}\mathrm{d}\vec{y} 
 \end{aligned}
 \end{equation}
 Hence due to \eqref{eq:grad-split-compactness_proof_2} and \eqref{eq:grad-split-compactness_proof_0_2} we have for    $\vec{\varphi}\in{\mathcal{D}}(\O;C^\infty_\sharp(Y^m))^m$, such that  $\left({\bf R} {\bf R}^T \nabla_y \right)\cdot \vec{\varphi}(\vec{x},\vec{y}) = 0$
 \begin{equation}
    \IOAY \left(\hat{\vec{\chi}}_0(\vec{x},\vec{y}) -  {\bf R} \nabla  u(\vec{x}) \right)
  \cdot\vec{\varphi}\left(\vec{x},\vec{y}\right)   \,\mathrm{d}\vec{x}\mathrm{d}\vec{y}  = 0
 \end{equation}
 We deduce, due to orthogonality in the dual pairing sense \eqref{eq:Green_identity}, that there exists  
 $u_1 \in L^p(\Omega;X_p)$ 
 such that 
 \begin{equation*}
     \hat{\vec{\chi}}_0(\vec{x},\vec{y}) =  {\bf R} \nabla  u(\vec{x}) +   {\bf R} {\bf R}^T\nabla_y  u_1(\vec{x},\vec{y})
 \end{equation*}
 We conclude that the limit of the gradient in \eqref{eq:grad-split-compactness_proof} becomes  
   \begin{equation}\label{eq:grad-split-compactness_proof_0_3}
  \begin{aligned}
 &
 \tilde{\chi}_0(\vec{x},\vec{y} =  {\bf R}^T   \hat{\chi}_0(\vec{x},\vec{y}) = {\bf R}^T   \left({\bf R} \nabla  u(\vec{x})  +  {\bf R} {\bf R}^T\nabla_y  u_1\left(\vec{x},\vec{y}\right) \right) = 
 \nabla  u(\vec{x})   + {\bf R}^T \nabla_y u_1\left(\vec{x},\vec{y}\right)  = \\ 
 &
 \mathrm{grad} \; u (\vec{x}) + \mathrm{grad}_{\bf R} \;  u_1\xy 
 \end{aligned}
 \end{equation}
 which completes the proof. \hfill \boxed{}

\vspace{5mm}
We define a strictly  monotone operator  $a$, which  satisfies the following assumptions,  i)-iv):
\begin{itemize}
    \item[i)] $ a(\cdot) $ is continuous  on $\R^n$
    \item[ii)] $ 0 \leq c_1  |\vec{\xi}|^p \leq  \left( a(\vec{\xi}) , \vec{\xi} \right) , \;  c_1>0, \; \forall \vec{\xi} \in \R^n$  
    
    \item[iii)]  $\left(a( \vec{\xi}_1)- a(\vec{\xi}_2), \vec{\xi}_1-\vec{\xi}_2 \right)  > 0   $ for all $ \vec{\xi}_1,\vec{\xi}_2  \in \R^n$, and almost every $\vec{x}\in \Omega$  and $\vec{y}\in\R^m$.
    \item[iv)]  $|a(\vec{\xi}) | \leq c_2 \left(1 + |\vec{\xi}|^{p-1}\right), \; c_2>0, \; \forall \vec{\xi} \in \R^n$
   
\end{itemize}

We  will have use of the following Lemma when  characterizing  the two-scale limit of divergences.

\begin{lemma}\label{lem:R-poisson_eqn}
Let $ 1 < p <\infty$, $1/p + 1/q =1$ and 
assume the  operator  $a$ satisfies assumptions i)-iv) above and that $f(\vec{x}, \cdot) \in L^q_{\sharp}(Y^m)$. The equation
\begin{equation}\label{eq:projected_local_eqn}
  -  \mathrm{div}_{\bf R} \; a\left( \mathrm{grad}_{\bf R}\; \theta( \vec{x}, \cdot) \right) = f(\vec{x}, \cdot) \;,  \qquad a.e. \;  \vec{x} \in \Omega
\end{equation}
with periodic boundary conditions,  has a unique weak solution  $\mathrm{grad}_{\bf R}\;{\theta}(\vec{x}, \cdot)$   in $L^{p}_{\sharp}(Y^m;\R^n)$,   
\end{lemma} 
{\bf Proof }
The proof follows from Browder (1963) and Minty (1963), \eg see \cite{lukkassen2002two}, page 62.
\hfill \boxed{}

\begin{proposition}\label{prop:div}
Let  $\{\vec{u}_\e\}$ be a uniformly bounded
sequence in ${\mathcal{W}}(\mathrm{div},\Omega)$.
Then there exist a subsequence 
 $\{\vec{u}_{\e_k}\}$ 
and functions in  $\vec{u}_0\in {\mathcal{W}^p}(\mathrm{div},\Omega; L^p( Y^m))$ 
and $\vec{u}_1\in L^p(\O, {\mathcal{W}}^p_{\sharp}(\mathrm{div}_{{\bf R}}, Y^m))$ 
such that
\begin{equation}\label{eq:divsplit_2s_weak}
\vec{u}_{\e_k} \wdto \vec{u}_0(\vec{x},\vec{y}), \qquad \mathrm{div}\; \vec{u}_{\e_k} \wdto  \mathrm{div} \; \vec{u}(\vec{x}) + \mathrm{div}_{\bf R} \; \vec{u}_1\xy,  \qquad {\e_k} \to 0
\end{equation}
with 
\begin{equation}\label{eqn:divfree_loocal}
 \mathrm{div}_{\bf R} \; \vec{u}_0\xy = 0
\end{equation}
and 
\begin{equation*}
\vec{u}(\vec{x}) = \int_{Y^m}  \vec{u}_0(\vec{x},\vec{y}) \; \mathrm{d}\vec{y}
\end{equation*}
$\vec{u}\in {\mathcal{W}^p}(\mathrm{div},\Omega)$. 
\end{proposition}
{\bf Proof }
The proof follows the lines of Lemma~5 and Proposition~6 in \cite{Wellander+etal2019} with appropriately  changed  function spaces.
 Let  $\phi\in L^q(\O)$  and $\psi\in L^q(\O,C_\sharp(Y^m))$. We have the weak limit of the divergence, 
 \begin{equation}\label{eq:div-split-compactness_proof_1}
  \begin{aligned}
 &
 \lim _{\e\to 0} \IOm    \nabla \cdot \vec{u}_{\e}(\vec{x})   \phi\left(\vec{x}\right) \,\mathrm{d}\vec{x} =   
\IOAY  \nabla \cdot \vec{u}_0(\vec{x},\vec{y})   \phi(\vec{x})\,\mathrm{d}\vec{x}\mathrm{d}\vec{y}    =   
\int_{\Omega}  \nabla \cdot \vec{u}(\vec{x})    \phi(\vec{x})\,\mathrm{d}\vec{x} \, ,\; \forall \phi\in L^q(\O)
 \end{aligned}
 \end{equation}
 where $\vec{u}(\vec{x})  = \int_{Y^m} \vec{u}_0(\vec{x},\vec{y}) \mathrm{d}\vec{y}$, where  $\vec{u}_0(\vec{x},\vec{y})$ is the two-scale cut-and-project limit with respect to ${\bf R}$. 
Next,  we have the corresponding  two-scale cut-and-project limit of the divergence   
  \begin{equation}\label{eq:div-split-compactness_proof_2}
  \begin{aligned}
 &
 \lim _{\e\to 0} \IOm    \nabla \cdot \vec{u}_{\e}(\vec{x})  \psi\left(\vec{x},\frac{ {\bf R}\vec{x}}{\e}\right) \,\mathrm{d}\vec{x} =  
\IOAY   {\chi}_0(\vec{x},\vec{y})   \psi(\vec{x},\vec{y})\,\mathrm{d}\vec{x}\mathrm{d}\vec{y}     \, ,\; \forall  \psi\in L^q(\O,C_\sharp(Y^m))
\end{aligned}
 \end{equation}
 It follows, after an integration by parts (twice) that 
  \begin{equation}\label{eq:div-free_limit-compactness_proof}
  \begin{aligned}
 &
0=  \lim _{\e\to 0} \e \IOm    \nabla \cdot \vec{u}_{\e}(\vec{x})  \psi\left(\vec{x},\frac{ {\bf R}\vec{x}}{\e}\right) \,\mathrm{d}\vec{x} =  
\IOAY  \mathrm{div}_{\bf R} \;   {\vec{u}}_0(\vec{x},\vec{y})   \psi(\vec{x},\vec{y})\,\mathrm{d}\vec{x}\mathrm{d}\vec{y}        \, ,\; \forall  \psi\in L^q(\O,C_\sharp(Y^m))
\end{aligned}
 \end{equation}
 which proves \eqref{eqn:divfree_loocal}. 
Define a function as the difference of the two-scale and the weak limits, \ie   $f(\vec{x},\vec{y}) := {\chi}_0(\vec{x},\vec{y})-   \nabla \cdot \vec{u}(\vec{x})$. We have $f(\vec{x},\cdot) \in L^p_{\sharp}(Y^m)$. 
Lemma~\ref{lem:R-poisson_eqn} implies that there is a unique
   $\mathrm{grad}_{\bf R}\;{\theta}(\vec{x}, \cdot)$  in $L^{q}_{\sharp}(Y^m;\R^n)$ that solves \eqref{eq:projected_local_eqn}, \ie  $ {\theta}(\vec{x}, \cdot) \in X_q$, defined in \eqref{eq:X-proj_space}. 
   Next, define  $\vec{u}_1(\vec{x},\vec{y}) := -    a\left(\mathrm{grad}_{\bf R}\;{\theta}(\vec{x}, \vec{y})\right)$ with $a$ as in Lemma~\ref{lem:R-poisson_eqn}. We get  
\begin{equation*}
    {\chi}_0(\vec{x},\vec{y}) =  \mathrm{div} \; \vec{u}(\vec{x}) + f(\vec{x},\vec{y}) =  \mathrm{div} \; \vec{u}(\vec{x})  + \mathrm{div}_{\bf R} u_1(\vec{x},\vec{y})   \; \in  L^p(\Omega \times Y^m) 
\end{equation*} 
which completes the proof.
 \hfill \boxed{}

\section{Homogenization of a quasiperiodic heterogeneous nonlinear electrostatic problem}

Let us now consider the quasiperiodic heterogeneous nonlinear electrostatic problem \eqref{eq:NL_static}. 
Standard estimates yield  solutions that are uniformly bounded in $W^{1,p}_0(\Omega)$ with respect to $\e$. We can now state the main homogenization result.
\begin{theorem}\label{theorem}
Let $\{u_\e\}$ be a sequence of solutions to   \eqref{eq:NL_static}.  The whole sequence $\{u_\e\}$  converges weakly in $W^{1,p}_0(\Omega)$  to the solution $\{u\}$ of the homogenized equation
   \begin{equation}  \label{eq:NL_static_hom}
     \left\{
        \begin{aligned}
          -  \mathrm{div} \; \int_{Y^m} \sigma\left( \vec{x}, \vec{y}, \nabla u(\vec{x}) + {\bf R}^T\nabla_{\vec{y}} u_1(\vec{x}, \vec{y}) \right)  \; \mathrm{d}\vec{y} & = f (\vec{x}) \;, \qquad \vec{x} \in \Omega\\
           \left. u \right|_{\partial\Omega} & = 0
        \end{aligned}
    \right.
 \end{equation}
 where ${\bf R}^T\nabla_{\vec{y}} u_1$ ($u_1 \in X_p$) is the unique solution of the local equation
 \begin{equation}  \label{eq:NL_static_cell}
          -  \mathrm{div}_{\vec{R}} \;  \sigma\left( \vec{x}, \vec{y}, \nabla u(\vec{x}) + {\bf R}^T\nabla_{\vec{y}} u_1(\vec{x}, \vec{y}) \right)  = 0   \qquad \vec{x} \in \Omega
 \end{equation}
\end{theorem}
{\bf Proof }
From the a priori estimates of sequences $u_\eta$ and $\sigma_\eta :=\sigma\left( \vec{x}, \frac{{\bf R}\vec{x}}{\e}, \nabla u_\e(\vec{x}) \right)$, there is subsequence such that $u_{\e_k} \dto u(\vec{x})$, $\mathrm{grad}\; u_{\e_k} \wdto \mathrm{grad} \; u (\vec{x}) + \mathrm{grad}_{\bf R} \;  u_1\xy$ and $\sigma_\eta\wdto\sigma_0(\vec{x},\vec{y})$,
when ${\e_k} \to 0$.
Since $f+\mathrm{div}\sigma_\eta=0$, equation  \eqref{eqn:divfree_loocal} yields   $\mathrm{div}_{{\bf R}}\sigma_0(\vec{x},\vec{y})=0$ and $f(\vec{x})+\mathrm{div}_{\vec{x}}\int_{Y^m}\sigma_0(\vec{x},\vec{y})\mathrm{d}\vec{y}=0$.  We now need to obtain an explicit
expression for $\sigma_0(\vec{x},\vec{y})$ in terms of $\sigma$, $u$ and $u_1$.
Following, \eg  \cite{Allaire1992} we introduce a test function $\psi_\eta(\vec{x})=\nabla\{u(\vec{x})+\eta\phi_1(\vec{x}, \frac{{\bf R}\vec{x}}{\e})\}+t\phi(\vec{x}, \frac{{\bf R}\vec{x}}{\e})$ where $t>0$, $\phi$ and $\phi_1$ are admissible test functions.
This ensures that $\psi_{\e_k} \wdto \mathrm{grad} \; u (\vec{x}) + \mathrm{grad}_{\bf R} \;  \phi_1\xy + t \phi\xy$.
Since $\sigma$ is strictly monotone, we have  
\begin{equation}
\int_\Omega\left\{\sigma_\eta-\sigma\left( \vec{x}, \frac{{\bf R}\vec{x}}{\e}, \psi_\eta\right)\right\}\cdot(\nabla u_\eta-\psi_\eta) \; \mathrm{d}\vec{x} \geq 0
\end{equation}
\ie
\begin{equation}
\int_\Omega\left\{-\mathrm{div}\sigma_\eta u_\eta-\sigma\left( \vec{x}, \frac{{\bf R}\vec{x}}{\e}, \psi_\eta\right)\cdot\nabla u_\eta -
\sigma_\eta \psi_\eta+\sigma\left( \vec{x}, \frac{{\bf R}\vec{x}}{\e}, \psi_\eta\right)\psi_\eta\right\} \; \mathrm{d}\vec{x} \geq 0
\end{equation}
Using (\ref{eq:NL_static}),  passing to the two-scale limit and using  the strong limit to get $u$ yields 
\begin{equation}
\begin{aligned}
&\int_\Omega\int_{Y^m}\left\{f(\vec{x})u(\vec{x})-\sigma\left( \vec{x}, \vec{y}, \psi_0(\vec{x},\vec{y})\right)\cdot(\mathrm{grad} \; u (\vec{x}) + \mathrm{grad}_{\bf R} \;  u_1\xy)-\sigma_0\xy\psi_0\xy\right. \nonumber \\
&+\left.\sigma\left( \vec{x}, \vec{y}, \psi_0(\vec{x},\vec{y})\right)\psi_0(\vec{x},\vec{y})\right\} \; \mathrm{d}\vec{x}\mathrm{d}\vec{y} \geq 0
\end{aligned}
\end{equation}
This equals, after a few integration by parts,   
\begin{equation}
\begin{aligned}
&\int_\Omega\int_{Y^m}\left\{f(\vec{x})u(\vec{x}) +  \mathrm{div} \; \sigma\left( \vec{x}, \vec{y}, \psi_0(\vec{x},\vec{y})\right) u (\vec{x}) - \sigma\left( \vec{x}, \vec{y}, \psi_0(\vec{x},\vec{y})\right) 
\cdot \mathrm{grad}_{\bf R} \;  u_1\xy \right. + \\
&\mathrm{div} \; \sigma_0\xy   u(\vec{x})   + \mathrm{div}_{\bf R}\sigma_0\xy \phi_1\xy  - \sigma_0\xy t \phi\xy - \\
& - \mathrm{div} \; \sigma\left( \vec{x}, \vec{y}, \psi_0(\vec{x},\vec{y})\right)  u(\vec{x}) -
\mathrm{div}_{\bf R} \; \sigma\left( \vec{x}, \vec{y}, \psi_0(\vec{x},\vec{y})\right)   \phi_1\xy +\\
& 
\left. 
\sigma\left( \vec{x}, \vec{y}, \psi_0(\vec{x},\vec{y})\right)  t \phi\xy   \right\} \; \mathrm{d}\vec{x}\mathrm{d}\vec{y} \geq 0
\end{aligned}
\end{equation}
The first terms in the first two rows cancel each other due to the statements above. We also note that the middle term in the second row vanishes.   The second term in the first row is canceled by the first term in the third row. 
Taking a sequence of functions $\mathrm{grad}_{\bf R} \;  \phi_1$ that converges strongly to $\mathrm{grad}_{\bf R} \;  u_1 $ in $L^p(\O, L^p(Y^n)^n))$, 
yields 
\begin{equation}
\begin{aligned}
&\int_\Omega\int_{Y^m}\left\{  
    \sigma (\vec{x}, \vec{y},\mathrm{grad} \; u (\vec{x}) + \mathrm{grad}_{\bf R} \;  u_1\xy+t\phi\xy) \cdot \mathrm{grad}_{\bf R} \;  u_1\xy  \right. + \\
&  - \sigma_0\xy t \phi\xy  -  
  \sigma (\vec{x}, \vec{y},\mathrm{grad} \; u (\vec{x}) + \mathrm{grad}_{\bf R} \;  u_1\xy+t\phi\xy) \cdot \mathrm{grad}_{\bf R} \;  u_1\xy  +\\
& 
\left. 
\sigma (\vec{x}, \vec{y},\mathrm{grad} \; u (\vec{x}) + \mathrm{grad}_{\bf R} \;  u_1\xy+t\phi\xy)  t \phi\xy   \right\} \; \mathrm{d}\vec{x}\mathrm{d}\vec{y} \geq 0
\end{aligned}
\end{equation}
The first row cancels the second term in the second row. We divide the two terms left by    $t>0$ and send   $t$   to zero and obtain 
 
\begin{equation}
\int_\Omega\int_{Y^m} \left[\sigma (\vec{x}, \vec{y},\mathrm{grad} \; u (\vec{x}) + \mathrm{grad}_{\bf R} \;  u_1\xy)-\sigma_0\xy\right]\phi\xy\mathrm{d}\vec{x}\mathrm{d}\vec{y} \geq 0
\end{equation} 
for all admissible test function, \eg   $\phi\in {\mathcal{D}}(\O;C^\infty_\sharp(Y^m))$.
It follows that $\sigma_0\xy=\sigma (\vec{x}, \vec{y},\mathrm{grad} \; u (\vec{x}) + \mathrm{grad}_{\bf R} \;  u_1\xy)$. Uniqueness of the solution of the limit equation (see e.g. \cite{lions1969quelques,wellander1998homogenization}) implies that the whole sequence converges. \hfill \boxed{}

\begin{proposition}[Correctors]\label{proposition}
If we assume that $u_1(\vec{x}, \vec{y})$ is smooth and $\sigma$ is uniformly monotone, then
\begin{equation*}
 \lim_{\eta \to 0} \left\|  \nabla  u_\e(\vec{x}) - \nabla \left\{u(\vec{x})+\eta u_1\left(\vec{x}, \frac{\vec{R}\vec{x}}{\e}\right)\right\}  \right\|_{L^p(\Omega;\R^n)} = 0
\end{equation*}
\end{proposition}
{\bf Proof } 
Considering $\psi_\eta(\vec{x})=\nabla\left\{u\left(\vec{x}\right)+\eta u_1\left(\vec{x}, \frac{\vec{R}\vec{x}}{\e}\right)\right\}$ and using that $\sigma$ is uniformly monotone yields
\begin{equation}
\int_\Omega\left\{\sigma_\eta-\sigma\left( \vec{x}, \frac{\vec{R}\vec{x}}{\e}, \psi_\eta\right)\right\}\cdot(\nabla u_\eta-\psi_\eta) \; \mathrm{d}\vec{x} \geq c
\int_\Omega { \mid\nabla u_\eta-\psi_\eta \mid}^p \; \mathrm{d}\vec{x}
\label{lastone}
\end{equation}
If follows from the fact that $\psi_\eta$ are admissible test functions that the left hand side of (\ref{lastone}) goes to zero in
$L^p(\O;\R^n)$. \hfill \boxed{}

\section{Concluding remarks}
We have applied two-scale cut-and-projection convergence to a canonical nonlinear electrostatic problem for quasiperiodic structures generated by a periodic geometry in a higher dimensional space. Compared with earlier work on homogenization of almost periodic monotone operators \cite{Braides1991,Nguetseng2007}, our annex problem has a simpler, less abstract  structure, and should therefore facilitate its numerical implementation in a variety of problems of physical interest, such as in electromagnetism \cite{wellander1998homogenization}, where intriguing features have been observed, such as transmitted femtosecond pulses developed a trailing diffusive exponential tail that led to some controversy \cite{ledermann2009multiple}.
We further note that our study can be adapted to the nonlinear elasticity case \cite{ponte1989overall}, whereby $\sigma$ would denote a rank-4 elasticity tensor.

\section*{Acknowledgement}

The work of E.C. was supported by  NSF grant DMS-1715680. S.G. acknowledges funding from the UK Engineering and
Physical Sciences Research Council (Grant No. EP/T002654/1)     

\bibliographystyle{cas-model2-names}


\bibliography{samplebib}



\end{document}